\documentclass[english,11pt]{article}
\usepackage{stmaryrd}
\usepackage{tipa}
\usepackage{tikz}
\usepackage{amsmath}
\usepackage{amssymb}
\usepackage{pst-node}
\usepackage{color}
\usepackage{extarrows}

\newtheorem{theorem}{Theorem}[section]
\newtheorem{corollary}[theorem]{Corollary}
\newtheorem{definition}[theorem]{Definition}
\newtheorem{conjecture}[theorem]{Conjecture}

\newtheorem{lemma}[theorem]{Lemma}

\newtheorem{proposition}[theorem]{Proposition}

\newtheorem{innercustomgeneric}{\customgenericname}
\providecommand{\customgenericname}{}
\newcommand{\newcustomtheorem}[2]{%
  \newenvironment{#1}[1]
  {%
   \renewcommand\customgenericname{#2}%
   \renewcommand\theinnercustomgeneric{##1}%
   \innercustomgeneric
  }
  {\endinnercustomgeneric}
}

\newcustomtheorem{customtheorem}{Theorem}
\newcustomtheorem{customlemma}{Lemma}

\newenvironment{proof}{\noindent {\bf Proof.}}{\rule{3mm}{3mm}\par\medskip}

\newcommand{\qed}{\fbox{\rule[1.2mm]{0.8mm}{0mm}}}

\newcommand{\Z}{\mbox{$\mathbb Z$}}

\newcommand{\setZTHR}{\mbox{$\langle\Z_3\rangle$}}

\begin{document}
\title{ Nowhere-zero $3$-flow of graphs  with small independence number}

 \author{ Jiaao Li\thanks{Department of Mathematics, West Virginia University, Morgantown, WV. 26506 USA.
{\em Email: joli@mix.wvu.edu}}
 , Rong Luo\thanks{Corresponding author, 
 School of Mathematics and Statistics,
Jiangsu Normal University,
  Xuzhou, Jiangsu, 221116,  China, 
  and 
  Department of Mathematics,
 West Virginia Univiersity,
 Morgantown, WV 26506 USA.
 {\em Email:  rongluo2007@gmail.com}}
 , Yi Wang\thanks{
 School of Mathematical Sciences,
Anhui University,
Hefei, Anhui 230601, China.
{\em Email: wangy@ahu.edu.cn}}
}

 \date{}
\maketitle

\begin{abstract}
Tutte's $3$-flow conjecture states that every $4$-edge-connected graph admits a nowhere-zero $3$-flow.  In this paper, we characterize all graphs with independence number at most $4$  that admit a nowhere-zero $3$-flow. The characterization of $3$-flow verifies Tutte's $3$-flow
conjecture for graphs with independence number at most $4$ and with order at least $21$. In addition, we prove that every odd-$5$-edge-connected graph with independence number at most $3$ admits a nowhere-zero $3$-flow.  To obtain these results, we introduce a new reduction method to handle odd wheels.
\end{abstract}

\section{Introduction}
Graphs in this paper are finite and loopless, but may contain parallel edges. We follow \cite{BoMu76} for undefined terms and notation.
For a graph $G$,  let $\alpha(G)$, $\kappa'(G)$,  and
$\delta(G)$  denote the independence number, the edge-connectivity,
and the minimum degree of $G$, respectively.
For vertex subsets $U,  W\subseteq V(G)$, let
$[U, W]_G = \{uw \in E(G)|u \in U, w \in W\}$. When $U =  \{u\}$ or
$W =\{w\}$, we use $[u,W]_G$ or $[U,w]_G$ for $[U,W]_G$, respectively. We also use $\partial_G(S)=[S,V(G)-S]_G$ to denote an edge-cut of $G$.
The subscript $G$ may be omitted when $G$ is understood from the context.

Let $D=D(G)$ be an orientation of $G$.
For each $v\in V(G)$, let $E^+_D(v)$ ($E^-_D(v)$, respectively) be the set of all arcs
directed out from (into, respectively) $v$. An {\em integer flow $(D, f)$} of $G$ is an orientation $D$ and a mapping $f: E(G) \mapsto \Z$  such that, for every vertex $v\in V(G)$,
\[
\sum_{e \in E^+_D(v)} f(e) - \sum_{e \in E^-_D(v)} f(e) =0.
\]
An integer flow $(D, f)$ is called a  {\em nowhere-zero $k$-flow} if $1\le|f(e)|\le k-1$,  for each edge $e\in E(G)$.

Let $d^+_D(v)=|E^+_D(v)|$ and $d^-_D(v)=|E^-_D(v)|$ denote
the out-degree and the in-degree of $v$ under the orientation $D$,
respectively.
A graph  $G$ admits a {\em modulo $3$-orientation},
or a mod $3$-orientation for short  if it has an orientation $D$ such that
$d^+_D(v) - d^-_D(v) \equiv 0\pmod3$ for every vertex $v \in V(G)$. It is well-known that a graph admits a nowhere-zero $3$-flow if and only if it admits a modulo $3$-orientation (see \cite{Tutte1949, Younger1983, Jaeg88}). Therefore, in this paper, we will study nowhere-zero $3$-flow in terms of modulo $3$-orientation.  The odd-edge-connectivity of a graph is  defined as the minimum size of an edge-cut of odd size. A  graph with low edge-connectivity may have high odd-edge-connectivity.

Tutte  posed the following famous $3$-Flow Conjecture, which appeared  in 1970s (see \cite{BoMu76}).

\begin{conjecture}\label{CON: 3flow}
  Every $4$-edge-connected graph admits a nowhere-zero $3$-flow.
\end{conjecture}

Thomassen \cite{Thomassen2012} settled the weak version of $3$-Flow Conjecture with edge-connectivity $8$ replacing $4$ and his result was further improved by Lov\'{a}sz, Thomassen, Wu and Zhang \cite{LTWZ13}.
\begin{theorem}\label{LTWZ}
{\em(Lov\'{a}sz et al. \cite{LTWZ13})} Every odd-$7$-edge-connected graph admits a nowhere-zero $3$-flow.
\end{theorem}

Jaeger et al. \cite{JLPT92} introduced the concept of group connectivity as generalizations of nowhere-zero flows.
 Let $Z(G,\Z_{3}) =\{b: V(G)\rightarrow \Z_{3}\;|$
$\sum_{v\in V(G)}b(v) \equiv 0  \pmod {3}\}$.
A graph $G$ is {\em   $\Z_{3}$-connected} if, for any
 $b\in Z(G,\Z_{3})$, there is an orientation $D$  such
that $d^+_D(v)-d^-_D(v)\equiv b(v)  \pmod {3}$ for every vertex $v \in V (G)$. Let  $\langle\Z_{3}\rangle$
denote   the family
of all $\Z_{3}$-connected graphs.  Jaeger et al. \cite{JLPT92} pointed out  that not every $4$-edge-connected graph is $\Z_{3}$-connected, and they further conjectured that
{\em  every $5$-edge-connected graph is $\Z_3$-connected.} This conjecture, if true,  implies Tutte's $3$-Flow Conjecture as Kochol \cite{Kochol01} showed that the $3$-Flow conjecture is equivalent to its restriction to $5$-edge-connected graphs.

 Luo et al. \cite{Luo2013} characterized graphs with independence number two that admit a  nowhere-zero $3$-flow.

\begin{theorem}\label{LuoMX2013}{\em (Luo et al. \cite{Luo2013})}
Let $G$ be a bridgeless graph with independence number $\alpha(G)\le 2$.
Then $G$ admits a nowhere-zero 3-flow if and only if G cannot be contracted to  $K_4$ or $G^3$,
and $G$ is not one of  three exceptional graphs, $G^3, G^{5}, G^{18}$ (see Figure \ref{FIG: 01}).
\end{theorem}
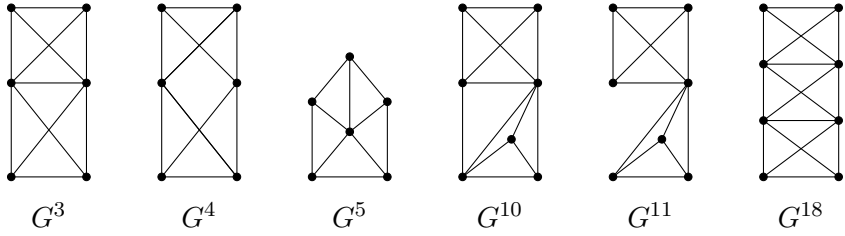
\begin{figure}[htp]
\begin{center}
\begin{tikzpicture}
\filldraw[black] (0,0) circle (0.05cm);
\filldraw[black] (0,1) circle (0.05cm);
\filldraw[black] (1,1) circle (0.05cm);
\filldraw[black] (1,0) circle (0.05cm);
\filldraw[black] (0,-1.25) circle (0.05cm);
\filldraw[black] (1,-1.25) circle (0.05cm);
\draw [-] (0,0)--(1,0);
\draw [-] (0,0)--(0,1);
\draw [-] (0,0)--(1,1);
\draw [-] (0,1)--(1,1);
\draw [-] (1,0)--(1,1);
\draw [-] (0,1)--(1,0);
\draw [-] (0,0)--(0,-1.25);
\draw [-] (0,0)--(1,-1.25);
\draw [-] (0,-1.25)--(1,0);
\draw [-] (1,0)--(1,-1.25);
\draw [-] (0,-1.25)--(1,-1.25);
\node at (0.5,-1.75){$G^3$};

\filldraw[black] (2,0) circle (0.05cm);
\filldraw[black] (2,1) circle (0.05cm);
\filldraw[black] (3,1) circle (0.05cm);
\filldraw[black] (3,0) circle (0.05cm);
\filldraw[black] (2,-1.25) circle (0.05cm);
\filldraw[black] (3,-1.25) circle (0.05cm);
\draw [-] (2,0)--(2,1);
\draw [-] (2,0)--(3,1);
\draw [-] (3,1)--(3,0);
\draw [-] (3,0)--(2,1);
\draw [-] (3,1)--(2,0);
\draw [-] (3,1)--(2,1);
\draw [-] (2,0)--(2,-1.25);
\draw [-] (2,0)--(3,-1.25);
\draw [-] (2,-1.25)--(3,0);
\draw [-] (2,0)--(3,-1.25);
\draw [-] (2,-1.25)--(3,-1.25);
\draw [-] (3,0)--(3,-1.25);
\node at (2.5,-1.75){$G^4$};

\filldraw[black] (4,-0.25) circle (0.05cm);
\filldraw[black] (4.5,0.35) circle (0.05cm);
\filldraw[black] (5,-0.25) circle (0.05cm);
\filldraw[black] (4,-1.25) circle (0.05cm);
\filldraw[black] (5,-1.25) circle (0.05cm);
\filldraw[black] (4.5,-0.65) circle (0.05cm);
\draw [-] (4,-0.25)--(4.5,0.35);
\draw [-] (4,-0.25)--(4,-1.25);
\draw [-] (5,-1.25)--(4,-1.25);
\draw [-] (5,-1.25)--(5,-0.25);
\draw [-] (5,-0.25)--(4.5,0.35);
\draw [-] (4.5,-0.65)--(4,-0.25);
\draw [-] (4.5,-0.65)--(4.5,0.35);
\draw [-] (4.5,-0.65)--(4,-1.25);
\draw [-] (4.5,-0.65)--(5,-0.25);
\draw [-] (4.5,-0.65)--(5,-1.25);
%\draw [-] (4,0)--(4,-1.25);
%\draw [-] (4,-1.25)--(5,0);
%\draw [-] (4,-1.25)--(5,-1.25);
%\draw [-] (4.65,-0.75)--(5,-1.25);
%\draw [-] (4,-1.25)--(4.65,-0.75);
%\draw [-] (5,0)--(4.65,-0.75);
%\draw [-] (5,0)--(5,-1.25);
\node at (4.5,-1.75){$G^{5}$};

\filldraw[black] (6,0) circle (0.05cm);
\filldraw[black] (6,1) circle (0.05cm);
\filldraw[black] (7,1) circle (0.05cm);
\filldraw[black] (7,0) circle (0.05cm);
\filldraw[black] (6,-1.25) circle (0.05cm);
\filldraw[black] (7,-1.25) circle (0.05cm);
\filldraw[black] (6.65,-0.75) circle (0.05cm);
\draw [-] (6,0)--(6,1);
\draw [-] (6,0)--(7,1);
\draw [-] (7,1)--(7,0);
\draw [-] (7,0)--(6,1);
\draw [-] (7,0)--(6,0);
\draw [-] (7,1)--(6,1);
\draw [-] (6,0)--(6,-1.25);
\draw [-] (6,-1.25)--(7,0);
\draw [-] (6,-1.25)--(7,-1.25);
\draw [-] (6.65,-0.75)--(7,-1.25);
\draw [-] (6,-1.25)--(6.65,-0.75);
\draw [-] (7,0)--(6.65,-0.75);
\draw [-] (7,0)--(7,-1.25);
\node at (6.5,-1.75){$G^{10}$};

\filldraw[black] (8,0) circle (0.05cm);
\filldraw[black] (8,1) circle (0.05cm);
\filldraw[black] (9,1) circle (0.05cm);
\filldraw[black] (9,0) circle (0.05cm);
\filldraw[black] (8,-1.25) circle (0.05cm);
\filldraw[black] (9,-1.25) circle (0.05cm);
\filldraw[black] (8.65,-0.75) circle (0.05cm);
\draw [-] (8,0)--(8,1);
\draw [-] (8,0)--(9,1);
\draw [-] (9,1)--(9,0);
\draw [-] (9,0)--(8,1);
\draw [-] (9,0)--(8,0);
\draw [-] (9,1)--(8,1);
\draw [-] (8,-1.25)--(9,0);
\draw [-] (8,-1.25)--(9,-1.25);
\draw [-] (8.65,-0.75)--(9,-1.25);
\draw [-] (8,-1.25)--(8.65,-0.75);
\draw [-] (9,0)--(8.65,-0.75);
\draw [-] (9,0)--(9,-1.25);
\node at (8.5,-1.75){$G^{11}$};

\filldraw[black] (10,0.25) circle (0.05cm);
\filldraw[black] (10,-0.5) circle (0.05cm);
\filldraw[black] (11,-0.5) circle (0.05cm);
\filldraw[black] (10,1) circle (0.05cm);
\filldraw[black] (11,1) circle (0.05cm);
\filldraw[black] (11,0.25) circle (0.05cm);
\filldraw[black] (10,-1.25) circle (0.05cm);
\filldraw[black] (11,-1.25) circle (0.05cm);
\draw [-] (10,1)--(10,-1.25);
\draw [-] (11,1)--(11,-1.25);
\draw [-] (10,1)--(11,1);
\draw [-] (10,-1.25)--(11,-1.25);
\draw [-] (10,1)--(11,0.25);
\draw [-] (11,1)--(10,0.25);
\draw [-] (10,0.25)--(11,0.25);
\draw [-] (10,0.25)--(11,-0.5);
\draw [-] (10,-0.5)--(11,-0.5);
\draw [-] (10,-0.5)--(11,.25);
\draw [-] (10,-1.25)--(11,-0.5);
\draw [-] (10,-0.5)--(11,-1.25);
\node at (10.5,-1.75){$G^{18}$};

\end{tikzpicture}
\caption{Graphs in Theorems \ref{LuoMX2013} and \ref{YangLi}}\label{FIG: 01}
\end{center}

\end{figure}

Yang et al. \cite{YLL16} further refined this result to characterize  $3$-edge-connected $\Z_3$-connected graphs with independence number two.
To state their theorem, we need to  introduce the concept of $\langle\Z_{3}\rangle$-reduction first.
Note that a  $K_1$ is $\Z_3$-connected, which is called a trivial $\Z_3$-connected graph, and thus for any graph $G$,
every vertex lies in a maximal $\Z_3$-connected subgraph of $G$.
Let $H_1,H_2, \ldots,H_c$ denote the collection of all
maximal $\Z_3$-connected subgraph of $G$.  We call $G' = G/(\cup^c_{i=1}E(H_i))$  the
{\em $\langle\Z_{3}\rangle$-reduction} of $G$, and  we say that $G$ is
{\em $\langle\Z_{3}\rangle$-reduced to $G'$}.
A graph $G$ is {\em $\langle\Z_{3}\rangle$-reduced}
if $G$ does not have a nontrivial $\Z_3$-connected subgraph.
By definition, the $\langle\Z_{3}\rangle$-reduction of a graph is always $\langle\Z_{3}\rangle$-reduced. It is shown in \cite{LaiH00} that a graph $G$ admits a nowhere-zero $3$-flow (is $\Z_3$-connected, respectively) if and only if its $\langle\Z_{3}\rangle$-reduction admits a nowhere-zero $3$-flow (is $\Z_3$-connected, respectively). Moreover, the potential minimal counterexamples to Conjecture \ref{CON: 3flow} must be $\langle\Z_{3}\rangle$-reduced graphs. Therefore in order to describe nowhere-zero $3$-flow and $\Z_3$-connectedness properties of certain family of graphs, it is sufficient to characterize all  $\langle\Z_{3}\rangle$-reductions of this family.

\begin{theorem}\label{YangLi}{\em (Yang et al. \cite{YLL16})}
  Let $G$ be a $3$-edge-connected graph with $\alpha(G)\le 2$. If $G$ is not one
of the $18$  graphs of order at most $8$, then $G$ is $\Z_3$-connected if and only if $G$ cannot be $\langle\Z_{3}\rangle$-reduced to one
of the graphs  in $\{K_4, G^3, G^4, G^{10},G^{11}\}$ (see Figure \ref{FIG: 01}).
\end{theorem}

 The purpose of this paper is to  further extend Theorem~\ref{LuoMX2013} to graphs with independence number at most $4$, and thus resolve the $3$-Flow Conjecture for  this family of graphs.

Denote  ${\cal F}_1 = \{H |$ $H$ is $\langle\Z_{3}\rangle$-reduced without mod $3$-orientation, $2\le |V(H)|\le 15$, $\alpha(H)\le 4$ and $\kappa'(H)\le 3\}$, and let ${\cal F}_2 = \{H |$ $H$ has no mod $3$-orientation and $14\le |V(H)|\le 20\}$.

\begin{theorem}\label{mod3} Let $G$ be a  graph  with $\alpha(G)\le 4$.
Then  $G$ admits a nowhere-zero $3$-flow if and only if
$G\notin {\cal F}_2$ and the $\langle\Z_{3}\rangle$-reduction of $G$ is not in ${\cal F}_1$.
\end{theorem}

 Since each graph in ${\cal F}_1$ is of edge-connectivity at most $3$, Theorem \ref{mod3} immediately leads the following, which verifies the $3$-Flow Conjecture for graphs with at least $21$ vertices and independence number at most $4$.

\begin{theorem}\label{mod3indle4}
 Every $4$-edge-connected graph $G$ with $|V(G)|\ge 21$ and $\alpha(G)\le 4$  admits a nowhere-zero $3$-flow.
\end{theorem}

 In Section 3, we will show that Theorem \ref{mod3} is equivalent to Theorem \ref{mod3indle4} (Lemma~\ref{reducetoorder15}).

For graphs with independence number at most $3$, we can  eliminate the order requirement in Theorem~\ref{mod3indle4} and prove the following theorem.

 \begin{theorem}\label{mod3indle3}
 Every $4$-edge-connected graph $G$ with  $\alpha(G)\le 3$  admits a nowhere-zero $3$-flow.
\end{theorem}

In fact, in Section 3, we will prove slightly stronger results than Theorems~\ref{mod3indle4} and ~\ref{mod3indle3} by replacing $4$-edge-connectivity with odd-$5$-edge-connectivity.

\begin{theorem}\label{mod3indle4odd}
 Every odd-$5$-edge-connected graph $G$ with $|V(G)|\ge 21$ and $\alpha(G)\le 4$  admits a nowhere-zero $3$-flow.
\end{theorem}

 \begin{theorem}\label{mod3indle3odd}
 Every odd-$5$-edge-connected graph $G$ with  $\alpha(G)\le 3$  admits a nowhere-zero $3$-flow.
\end{theorem}

{\bf Remark.} There are quite a few graphs in the family ${\cal F}_1$ that are far from being described by hand.  In particular, the $18$ special graphs of order at most $8$ demonstrated by Yang et al. \cite{YLL16} can be modified to construct graphs in ${\cal F}_1$ by replacing a vertex of $K_4$ with one of those graphs. Also, many graphs obtained from  $2$-sum of two small non-$3$-flow admissible graphs are in ${\cal F}_1$.

While some splitting technique can not be applied for $\Z_3$-connectedness, it seems  very complicated to obtain  analogous results for  $\Z_3$-connectedness of graphs with small independence number
 via modifying the method of this paper and  much more involved discussion on small graphs are needed. However, such characterization for $\Z_3$-connectedness is interesting.  Note that Jaeger et al. \cite{JLPT92} constructed
a $4$-edge-connected graph $G$ of order $12$ with $\alpha(G)=3$,
which is not $\Z_3$-connected.

The organization of the rest of the paper is as follows:
Tools and preliminaries will be given in  Section 2 and the proofs of the main results will be presented in Section 3.

\section{Preliminaries}
In this section, we  display and develop some  tools needed in the proofs of the main results.
\subsection{Tools}

Lemma \ref{kn} is a summary of certain basic properties from \cite{JLPT92, LaiH00}.
\begin{lemma} \label{kn} Let $G$ be a graph. Each of the following holds:
\\
(i) If $G\in\langle\Z_{3}\rangle$ and
$e\in E(G)$, then $G/e\in\langle\Z_{3}\rangle$.
\\
(ii) If $H\subseteq G$,  and if both $H\in\langle\Z_{3}\rangle$ and $G/H \in\langle\Z_{3}\rangle$, then $G\in\langle\Z_{3}\rangle$.
\\
(iii) $G$ admits a mod $3$-orientation if and only if its  $\langle\Z_{3}\rangle$-reduction does.
\\
(iv) $G\in \langle\Z_{3}\rangle$ if and only if its  $\langle\Z_{3}\rangle$-reduction is $K_1$.
\\
(v) A cycle $C_n$ is $\Z_3$-connected if and only if $n=2$.
\\
(vi) The complete graph $K_n$ is $\Z_3$-connected if and only if $n=1$ or $n \geq 5$.
\end{lemma}

\begin{lemma}\label{min5}{\em (Han et al.  \cite{HanLL16})}
  Every $\setZTHR$-reduced graph has minimum degree at most $5$.
\end{lemma}

It has been  extensively  studied on the graphs admitting
nowhere-zero $3$-flows or  being $\Z_{3}$-connected under degree  conditions.  For example, Barat and Thomassen \cite{BaTh06} presented some degree conditions to ensure a simple
graph to be $\mathbb{Z}_3$-connected.
Fan and Zhou \cite{FaZh08} and Luo et al. \cite{Luo2008} characterized graphs admitting nowhere-zero $3$-flow and all $\Z_3$-connected graphs under Ore-condition, respectively, where  a simple graph $G$ satisfies {\bf Ore-condition}, if  for every pair of nonadjacent vertices $u$ and $v$ in $G$, $d_G(u) + d_G(v) \ge |V(G)|$.  Their results will be needed in our proofs to handle small graphs.

\begin{figure}[htp]
\begin{center}
\begin{tikzpicture}
\filldraw[black] (0,0) circle (0.05cm);
\filldraw[black] (0,1) circle (0.05cm);
\filldraw[black] (1,1) circle (0.05cm);
\filldraw[black] (1,0) circle (0.05cm);
\filldraw[black] (0,-1.25) circle (0.05cm);
\filldraw[black] (1,-1.25) circle (0.05cm);
\draw [-] (0,0)--(1,0);
\draw [-] (0,0)--(0,1);
\draw [-] (0,0)--(1,1);
\draw [-] (0,1)--(1,1);
\draw [-] (0,1)--(1,-1.25);
\draw [-] (0,1)--(1,0);
\draw [-] (0,-1.25)--(1,1);
\draw [-] (0,0)--(1,-1.25);
\draw [-] (0,-1.25)--(1,0);
%\draw [-] (1,0)--(1,-1.25);
\draw [-] (0,-1.25)--(1,-1.25);
\node at (0.5,-1.75){$(1)$};

\filldraw[black] (2,0) circle (0.05cm);
\filldraw[black] (2,1) circle (0.05cm);
\filldraw[black] (3,1) circle (0.05cm);
\filldraw[black] (3,0) circle (0.05cm);
\filldraw[black] (2,-1.25) circle (0.05cm);
\filldraw[black] (3,-1.25) circle (0.05cm);
\draw [-] (2,0)--(2,1);
\draw [-] (2,0)--(3,1);
\draw [-] (2,0)--(3,0);
\draw [-] (3,1)--(3,0);
\draw [-] (3,0)--(2,1);
\draw [-] (3,1)--(2,0);
\draw [-] (3,1)--(2,1);
\draw [-] (2,0)--(2,-1.25);
\draw [-] (2,0)--(3,-1.25);
\draw [-] (2,-1.25)--(3,0);
\draw [-] (2,0)--(3,-1.25);
\draw [-] (2,-1.25)--(3,-1.25);
\draw [-] (3,0)--(3,-1.25);
\node at (2.5,-1.75){$(2)$};

\filldraw[black] (4,-0.25) circle (0.05cm);
\filldraw[black] (4.5,0.35) circle (0.05cm);
\filldraw[black] (5,-0.25) circle (0.05cm);
\filldraw[black] (4,-1.25) circle (0.05cm);
\filldraw[black] (5,-1.25) circle (0.05cm);
\filldraw[black] (4.5,-0.65) circle (0.05cm);
\draw [-] (4,-0.25)--(4.5,0.35);
\draw [-] (4,-0.25)--(4,-1.25);
\draw [-] (5,-1.25)--(4,-1.25);
\draw [-] (5,-1.25)--(5,-0.25);
\draw [-] (5,-0.25)--(4.5,0.35);
\draw [-] (4.5,-0.65)--(4,-0.25);
\draw [-] (4.5,-0.65)--(4.5,0.35);
\draw [-] (4.5,-0.65)--(4,-1.25);
\draw [-] (4.5,-0.65)--(5,-0.25);
\draw [-] (4.5,-0.65)--(5,-1.25);
%\draw [-] (4,0)--(4,-1.25);
%\draw [-] (4,-1.25)--(5,0);
%\draw [-] (4,-1.25)--(5,-1.25);
%\draw [-] (4.65,-0.75)--(5,-1.25);
%\draw [-] (4,-1.25)--(4.65,-0.75);
%\draw [-] (5,0)--(4.65,-0.75);
%\draw [-] (5,0)--(5,-1.25);
\node at (4.5,-1.75){$(3)$};

\filldraw[black] (6.5,-0.5) circle (0.05cm);
\filldraw[black] (6,1) circle (0.05cm);
\filldraw[black] (7,1) circle (0.05cm);
\filldraw[black] (6.5,0.25) circle (0.05cm);
\filldraw[black] (6,-1.25) circle (0.05cm);
\filldraw[black] (7,-1.25) circle (0.05cm);

\draw [-] (6,-1.25)--(6,1);
\draw [-] (7,-1.25)--(7,1);
\draw [-] (6,1)--(7,1);
\draw [-] (6,-1.25)--(7,-1.25);
\draw [-] (6.5,0.25)--(6,1);
\draw [-] (6.5,0.25)--(7,1);
\draw [-] (6.5,-0.5)--(6,-1.25);
\draw [-] (6.5,-0.5)--(7,-1.25);
\draw [-] (6.5,-0.5)--(6.5,0.25);

\node at (6.5,-1.75){$(4)$};

\filldraw[black] (8,0) circle (0.05cm);
\filldraw[black] (9,0) circle (0.05cm);
\filldraw[black] (8,-1.25) circle (0.05cm);
\filldraw[black] (9,-1.25) circle (0.05cm);

\draw [-] (8,0)--(9,0);
\draw [-] (8,0)--(8,-1.25);
\draw [-] (9,0)--(8,-1.25);
\draw [-] (8,0)--(9,-1.25);
\draw [-] (9,0)--(9,-1.25);
\draw [-] (8,0)--(8,-1.25);
\draw [-] (9,-1.25)--(8,-1.25);
\node at (8.5,-1.75){$(5)$};

\filldraw[black] (10,-0.25) circle (0.05cm);
\filldraw[black] (11,-0.25) circle (0.05cm);
\filldraw[black] (10,-1.25) circle (0.05cm);
\filldraw[black] (11,-1.25) circle (0.05cm);
\filldraw[black] (10.5,0.5) circle (0.05cm);

\draw [-] (10,-0.25)--(10,-1.25);
\draw [-] (11,-0.25)--(10,-1.25);
%\draw [-] (10,-0.25)--(11,-0.25);
\draw [-] (10,-1.25)--(11,-1.25);
\draw [-] (10,-0.25)--(11,-1.25);
\draw [-] (11,-1.25)--(11,-0.25);
\draw [-] (9,-1.25)--(8,-1.25);
\draw [-] (10.5,0.5)--(10,-0.25);
\draw [-] (10.5,0.5)--(11,-0.25);
\node at (10.5,-1.75){$(6)$};

%%%%%%%%%%%%%%%%%%%%%%%

\filldraw[black] (0,-4) circle (0.05cm);
\filldraw[black] (0,-3) circle (0.05cm);
\filldraw[black] (1,-3) circle (0.05cm);
\filldraw[black] (1,-4) circle (0.05cm);
\filldraw[black] (0,-5.25) circle (0.05cm);
\filldraw[black] (1,-5.25) circle (0.05cm);
\draw [-] (0,-4)--(1,-4);
%\draw [-] (0,-4)--(0,-3);
\draw [-] (0,-4)--(1,-3);
\draw [-] (0,-3)--(1,-3);
\draw [-] (0,-3)--(1,-5.25);
\draw [-] (0,-3)--(1,-4);
\draw [-] (0,-5.25)--(1,-3);
\draw [-] (0,-4)--(1,-5.25);
\draw [-] (0,-5.25)--(1,-4);
%\draw [-] (1,0)--(1,-1.25);
\draw [-] (0,-5.25)--(1,-5.25);
\node at (0.5,-5.75){$(7)$};

\filldraw[black] (2,-4) circle (0.05cm);
\filldraw[black] (2,-3) circle (0.05cm);
\filldraw[black] (3,-3) circle (0.05cm);
\filldraw[black] (3,-4) circle (0.05cm);
\filldraw[black] (2,-5.25) circle (0.05cm);
\filldraw[black] (3,-5.25) circle (0.05cm);
\draw [-] (2,-4)--(2,-3);
\draw [-] (2,-4)--(3,-3);
%\draw [-] (2,-4)--(3,-4);
\draw [-] (3,-3)--(3,-4);
\draw [-] (3,-4)--(2,-3);
\draw [-] (3,-3)--(2,-4);
\draw [-] (3,-3)--(2,-3);
\draw [-] (2,-4)--(2,-5.25);
\draw [-] (2,-4)--(3,-5.25);
\draw [-] (2,-5.25)--(3,-4);
\draw [-] (2,-4)--(3,-5.25);
\draw [-] (2,-5.25)--(3,-5.25);
\draw [-] (3,-4)--(3,-5.25);
\node at (2.5,-5.75){$(8)$};

\filldraw[black] (4,-4.25) circle (0.05cm);
\filldraw[black] (4.5,-3.65) circle (0.05cm);
\filldraw[black] (5,-4.25) circle (0.05cm);
\filldraw[black] (4,-5.25) circle (0.05cm);
\filldraw[black] (5,-5.25) circle (0.05cm);
%\filldraw[black] (4.5,-0.65) circle (0.05cm);
%\draw [-] (4,-0.25)--(4.5,0.35);
\draw [-] (4,-4.25)--(4,-5.25);
\draw [-] (5,-5.25)--(4,-5.25);
\draw [-] (5,-5.25)--(5,-4.25);
\draw [-] (5,-4.25)--(4.5,-3.65);
\draw [-] (4,-4.25)--(4.5,-3.65);
\draw [-] (4,-4.25)--(5,-5.25);
\draw [-] (5,-4.25)--(4,-5.25);
\draw [-] (4,-4.25)--(5,-4.25);

\node at (4.5,-5.75){$(9)$};

\filldraw[black] (6,-4) circle (0.05cm);
\filldraw[black] (7,-4) circle (0.05cm);
\filldraw[black] (6,-5.25) circle (0.05cm);
\filldraw[black] (7,-5.25) circle (0.05cm);

\draw [-] (6,-4)--(7,-4);
\draw [-] (6,-4)--(7,-5.25);
%\draw [-] (7,-4)--(6,-5.25);
%\draw [-] (6,-4)--(7,-5.25);
\draw [-] (7,-4)--(7,-5.25);
\draw [-] (6,-4)--(6,-5.25);
\draw [-] (7,-5.25)--(6,-5.25);
\node at (6.5,-5.75){$(10)$};

\filldraw[black] (8,-4) circle (0.05cm);
\filldraw[black] (9,-4) circle (0.05cm);
\filldraw[black] (8,-5.25) circle (0.05cm);
\filldraw[black] (9,-5.25) circle (0.05cm);

\draw [-] (8,-4)--(9,-4);
%\draw [-] (8,-4)--(8,-5.25);
%\draw [-] (9,-4)--(8,-5.25);
%\draw [-] (8,-4)--(9,-5.25);
\draw [-] (9,-4)--(9,-5.25);
\draw [-] (8,-4)--(8,-5.25);
\draw [-] (9,-5.25)--(8,-5.25);
\node at (8.5,-5.75){$(11)$};

\filldraw[black] (10,-5.25) circle (0.05cm);
\filldraw[black] (11,-5.25) circle (0.05cm);
\filldraw[black] (10.5,-4) circle (0.05cm);

\draw [-] (10,-5.25)--(11,-5.25);

\draw [-] (10.5,-4)--(11,-5.25);
\draw [-] (10.5,-4)--(10,-5.25);
\node at (10.5,-5.75){$(12)$};

\end{tikzpicture}
\caption{Graphs in Theorems \ref{fanzhou} and \ref{luoorez3}}\label{FIG: 02}
\end{center}

\end{figure}
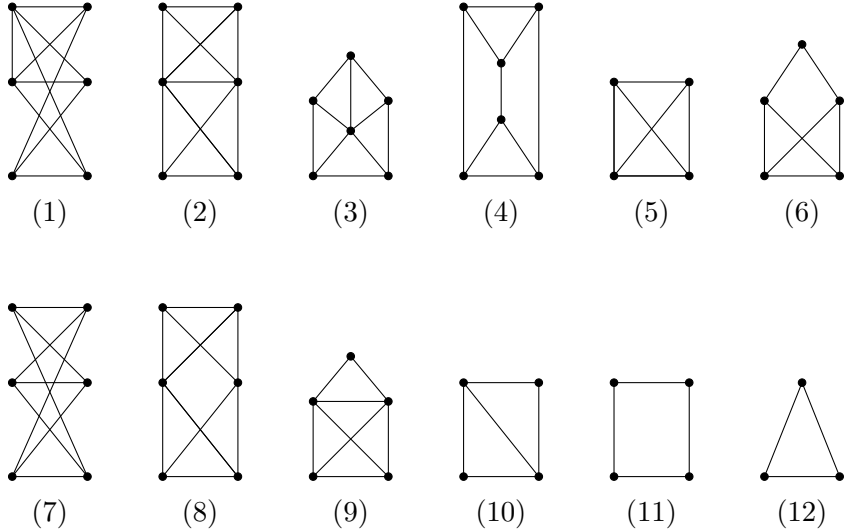

\begin{theorem}{\em(Fan and Zhou \cite{FaZh08})}
\label{fanzhou}
Let $G$ be a simple graph on $n \geq 3$ vertices satisfying the Ore-condition. Then $G$ admits a nowhere-zero $3$-flow except for six specified small graphs (see Figure \ref{FIG: 02} (1)-(6)).
\end{theorem}

\begin{theorem}\label{luoorez3}{\em(Luo et al. \cite{Luo2008})} Let $G$ be a simple graph on $n \geq 3$ vertices satisfying the Ore-condition. Then $G$ is $\mathbb{Z}_{3}$-connected except for $12$ specified small graphs (see Figure \ref{FIG: 02} (1)-(12)).
\end{theorem}

%%%%%%%%%%%%%%%%%%%%%%%%%%%%%%%%%%%%%%%%%%%%%%%%%%%%%%%%%%%%%%%%%%%%%%%%%%%%%%%%%%%%%%%%%%%%%%

 Let $u_1v$ and $u_2v$ be two distinct edges in $G$. Denote $G_{[v, u_1u_2]}$ to  be the graph obtained from $G$ by deleting the edges $u_1v, u_2v$ and adding a new edge $u_1u_2$, which is called the {\em lifting operation} (see \cite{Thomassen2012, LTWZ13}).  The following splitting lemma of Zhang \cite{CQsplitting02} shows that the odd-edge-connectivity is preserved under certain lifting operation.

\begin{lemma}{\em (Zhang \cite{CQsplitting02})}
\label{splitting}
  Let $G$ be a graph with odd-edge-connectivity $k$. Assume there is a vertex $v\in V(G)$ with $d(v)\neq k$ and $d(v)\neq 2$. Then there exists a pair of edges $u_1v, u_2v$ in $\partial_G(v)$ such that $G_{[v, u_1u_2]}$ preserves odd-edge-connectivity $k$.
\end{lemma}

\medskip
\noindent
{\bf Remark.} Lemma \ref{splitting}  does not apply to  a vertex $v$ of degree two as $v$ is an isolated vertex in $G_{[v, u_1u_2]}$. While in most of the flow problems (include the proofs in this paper), a degree two vertex $v$ does not appear in the minimal counterexamples since we could apply induction on the graph obtained from $G_{[v, u_1u_2]}$ by deleting the isolated vertex $v$. For this reason, we frequently ignore the discussion of degree two vertices when we  apply Lemma \ref{splitting}.

\subsection{A new contraction method to handle odd wheels}
A wheel $W_n$ is the graph obtained from an $n$-cycle by
adding a new vertex, called the center of the wheel, which is joined to every vertex of the $n$-cycle.
$W_n$ is odd (or even, respectively)  if $n$ is odd (or even, respectively). The complete graph $K_4$ can be viewed as a $W_3$.

\begin{lemma}\label{neq0b} Let $k$ be a positive integer.
\\
(i) {\em (DeVos et al. \cite{DeXY06})} Every even wheel $W_{2k}$ is $\Z_3$-connected.
\\
(ii) {\em (Xu \cite{XuRu04})} Let $b\in Z(W_{2k+1}, \Z_3)$. If there exists  $b(v)\neq 0$ for some $v\in V(W_{2k+1})$, then there is an orientation $D$ of $W_{2k+1}$ such that $d_D^+(x)-d_D^-(x)\equiv b(x)\pmod 3$ for any $x\in V(W_{2k+1})$.
\end{lemma}

\begin{figure}
\begin{center}
\begin{tikzpicture}
\filldraw[black] (0,5) circle (0.05cm);
\filldraw[black] (0,3.5) circle (0.05cm);
\filldraw[black] (1.5,5) circle (0.05cm);
\filldraw[black] (1.5,3.5) circle (0.05cm);
\draw [-] (0,5)--(0,3.5);
\draw [-] (0,5)--(1.5,5);
\draw [-] (0,5)--(1.5,3.5);
\draw [-] (0,3.5)--(1.5,3.5);
\draw [-] (1.5,5)--(1.5,3.5);
\draw [-] (1.5,5)--(0,3.5);
\draw [-] (0,5)--(-0.8,5.8);
\draw [-] (0,5)--(-0.5,6);
\draw [-] (0,3.5)--(-0.8,3.3);
\draw [-] (0,3.5)--(-0.5,3.5);
\draw [-] (1.5,3.5)--(2,3.5);
\draw [-] (1.5,3.5)--(2.3,3.3);
\draw [-] (1.5,3.5)--(2.6,3);
\draw [-] (1.5,5)--(2.15,5.9);
\node at (3,4.25){$\Longrightarrow$};
\filldraw[black] (5,4.25) circle (0.05cm);
\filldraw[black] (6.5,4.25) circle (0.05cm);
\draw [-] (5,4.25)--(4.2,4.95);
\draw [-] (5,4.25)--(3.9,4.75);
\draw [-] (5,4.25)--(4.2,3.45);
\draw [-] (5,4.25)--(4.5,3.25);
\draw [-] (5,4.25)--(6.5,4.25);
\draw [-] (5,4.25)--(4.2,4.95);
\draw [-] (6.5,4.25)--(7.15,5.15);
\draw [-] (6.5,4.25)--(7.6,3.75);
\draw [-] (6.5,4.25)--(7.3,3.45);
\draw [-] (6.5,4.25)--(7,3.25);
\draw[dashed] (0,4.25)
ellipse (0.4cm and 1.1cm);
\draw[dashed] (1.5,4.25)
ellipse (0.4cm and 1.1cm);

\filldraw[black] (0,0) circle (0.05cm);
\filldraw[black] (-1.2,-0.2) circle (0.05cm);
\filldraw[black] (1.2,-0.2) circle (0.05cm);
\filldraw[black] (0,-1.2) circle (0.05cm);
\filldraw[black] (-1,1) circle (0.05cm);
\filldraw[black] (1,1) circle (0.05cm);
\draw [-] (0,0)--(-1,1);\draw [-] (0,0)--(1,1);\draw [-] (0,0)--(-1.2,-0.2);\draw [-] (0,0)--(1.2,-0.2);\draw [-] (0,0)--(0,-1.2);
\draw [-] (1,1)--(-1,1);\draw [-] (-1.2,-0.2)--(-1,1);\draw [-] (0,-1.2)--(-1.2,-0.2);\draw [-] (0,-1.2)--(1.2,-0.2);\draw [-] (1,1)--(1.2,-0.2);
\draw [-] (-1.8,1.8)--(-1,1);
\draw [-] (1.8,1.8)--(1,1);
\draw [-] (2,1.6)--(1,1);\draw [-] (-1.2,-0.2)--(-1.8,-1);\draw [-] (-1.2,-0.2)--(-1.6,-1.2);\draw [-] (0,-1.2)--(-0.2,-2);\draw [-] (0,-1.2)--(0.2,-2);\draw [-] (1.2,-0.2)--(2,-1);
\node at (3.2,-0.2){$\Longrightarrow$};
\filldraw[black] (5.5,-0.2) circle (0.05cm);
\filldraw[black] (7,-0.2) circle (0.05cm);
\draw [-] (5.5,-0.2)--(7,-0.2);
\draw [-] (5.5,-0.2)--(4.7,0.6); \draw [-] (5.5,-0.2)--(4.7,-1); \draw [-] (5.5,-0.2)--(4.9,-1.2);
\draw [-] (7.8,0.6)--(7,-0.2); \draw [-] (8,0.4)--(7,-0.2); \draw [-] (8,-0.2)--(7,-0.2); \draw [-] (8,-0.8)--(7,-0.2); \draw [-] (7.8,-1)--(7,-0.2);
\draw[dashed] (-1.1,0.5)
ellipse (0.4cm and 1cm);
\draw[dashed] (0.55,-0.1)
ellipse (0.9cm and 1.7cm);
\end{tikzpicture}
\caption{Example of contraction in Lemma \ref{oddwheel} }\label{FIG: contraction}
\end{center}

\end{figure}
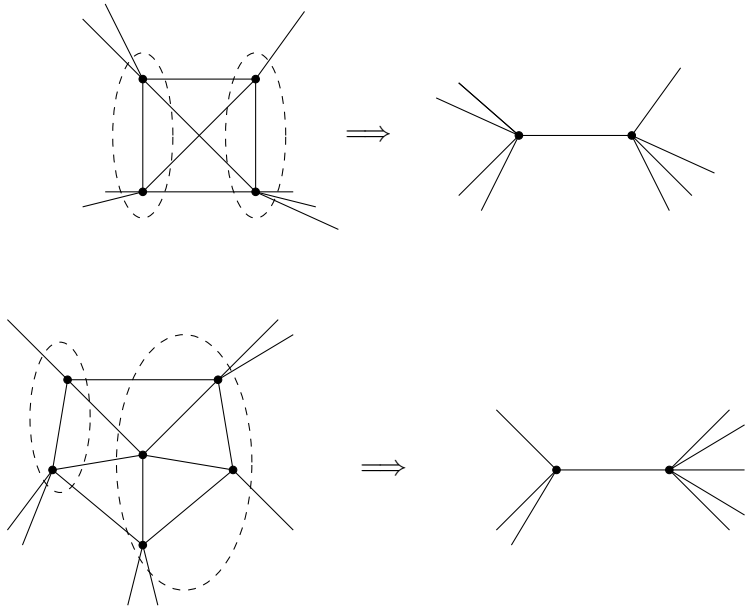

Lemma~\ref{neq0b}(ii) tells that an odd wheel is almost $\Z_3$-connected except when the boundary $b\in Z(W_{2k+1}, \Z_3)$  is a constant  zero function.
Thus  if a graph contains an odd wheel and if  the resulting graph   admits a nowhere-zero $3$-flow (is $\Z_3$-connected, respectively) after contracting an odd wheel into a $K_2$, then  so does (so is, respectively) the original graph.
Therefore we have the following lemma.
\begin{lemma}\label{oddwheel}
  Let $G$ be a  connected graph that contains a $W_{2k+1}$ as a proper subgraph of $G$.  Let $X, Y$ be a partition of $V(W_{2k+1})$, and let $G_{[X, Y]}$ be the graph obtained  from $G$ by deleting the edges of $E(W_{2k+1})$, contracting $X$ and $Y$ into $x$ and $y$, respectively, and adding a new edge $xy$ (see Figure \ref{FIG: contraction}).
\\
  (i) If  $G_{[X, Y]}$ has a mod $3$-orientation, then so does $G$.
\\
  (ii) If $G_{[X, Y]}$ is $\Z_3$-connected, then $G$ is $\Z_3$-connected.
\end{lemma}

An edge cut $\partial_G(S)=[S,V(G)-S]$ in a connected graph
$G$ is {\em essential} if at least two components of $G-\partial_G(S)$ are nontrivial, where a component is called {\em nontrivial} if it contains at least one edge. A graph is {\em essentially $k$-edge-connected} if it does not
have an essential edge cut with fewer than $k$ edges. Observe that, in a highly essentially edge-connected graph, if we contract an odd wheel into a single edge as described  in Lemma \ref{oddwheel}, then the edge connectivity of the  resulting graph cannot  drop too much. To formulate this for later application, we define the following special contraction of odd wheels.

\begin{definition}\label{DEF: Wcontraction}
 Let $G$ be a  connected graph and $ W_{2k+1}$ be a proper subgraph of $G$. $G_1=G_{[X, Y]}$ is a {\bf $W$-contraction} of $G$ if $X,Y$ form a partition of $V(W_{2k+1})$ and  one of $X, Y$ consists of two adjacent vertices in the $(2k+1)$-cycle of $W_{2k+1}$ (see Figure \ref{FIG: contraction}).
\end{definition}

Note  that, in a $W$-contraction of $G$, the original $4$ edges in $[X,Y]_G$  are replaced by a single edge $K_2=xy$. Hence an essential edge-cut of size $k$ in $G$ results  in an edge-cut of size at least $k-3$ in the $W$-contraction. It is also obvious that any $W$-contraction of $G$ has minimal degree at least $5$ provided that $G$ is $5$-edge-connected. Therefore, we obtain the following proposition.

\begin{proposition}\label{PROP: Wcontraction}
Let $G$ be  a $5$-edge-connected essentially $8$-edge-connected graph. If $G$ contains an odd wheel as a proper subgraph, then every $W$-contraction of $G$ remains $5$-edge-connected.
\end{proposition}

\subsection{Small graphs}

We shall discuss certain graphs of small order to serve for the induction basis of the proofs.

Denote $r(n, \Z_3) = \max\{|E(G)| ~|~ \mbox{$|V(G)| = n$ and $G$ is $\langle\Z_{3}\rangle$-reduced}\}$.  We determine $r(n, \Z_3)$ when $n$ is small in the following, which is needed in  later proof.
\begin{lemma} \label{reducedgraph} $r(1, \Z_3)=0$, $r(2, \Z_3)=1$, $r(3, \Z_3)=3$, $r(4, \Z_3)=6$, $r(5, \Z_3)=8$, $r(6, \Z_3)=11$, $r(7, \Z_3)=13$.
\end{lemma}
\proof~ Since a $\langle\Z_{3}\rangle$-reduced graph is simple by Lemma \ref{kn}(v),  it is routine to compute $r(n, \Z_3)$ when $n\le 4$. For $n=5$,  $K_5 -e$  is not $\langle\Z_{3}\rangle$-reduced  for any edge $e$ in $K_5$ because  it contains a $\Z_3$-connected subgraph, namely the wheel $W_4$ (by Lemma \ref{neq0b}(i)). Howerver, it is straightforward to show that $K_5$ deleting two incident edges is  $\langle\Z_{3}\rangle$-reduced (see Figure \ref{FIG: 02} (9)). Therefore $r(5, \Z_3)=8$.  We are to show $r(6, \Z_3)=11$ and $r(7, \Z_3)=13$ below.

 Let $G$ be a $\langle\Z_{3}\rangle$-reduced graph of order $6$. Since  every subgraph of a $\langle\Z_{3}\rangle$-reduced graph is also $\langle\Z_{3}\rangle$-reduced,  we have $|E(G)| \leq \delta(G) + r(5, \Z_3) = \delta(G) + 8$. By Lemma  \ref{kn}(v), $G$ is simple.  If $\delta(G) \leq 2$, then $|E(G)| \leq 10$. If $\delta(G) \geq 3$, then $G$ satisfies Ore-condition and thus by Theorem \ref{luoorez3}, we have $|E(G)|\le 11$ with equality  if and only if    $G$ is isomorphic to $G^3$ in Figure \ref{FIG: 01}.  This proves $r(6,\Z_3) = 11$ and  the only $\langle\Z_{3}\rangle$-reduced graph of order $6$ with  $11$ edges is $G^3$.

 Clearly, the graph obtained from $G^3$ by adding a new vertex with two nonparallel edges connecting to $G^3$ is $\langle\Z_{3}\rangle$-reduced. So $r(7, \Z_3)\ge 13$.  Let $G$ be a $\langle\Z_{3}\rangle$-reduced graph of order $7$. Then, by Lemma \ref{kn}(v), $G$ is  simple and by Theorem \ref{luoorez3},  $\delta(G)\le 3$. If $\delta(G)\le 2$, then $|E(G)|\le 2+ r(6, \Z_3)\le 13$. Assume $\delta(G)=3$.  Then $|E(G)| \leq 14$. If $|E(G)| = 14$, then $G-v = G^3$ for any degree $3$ vertex $v$.  Then $G$ must be the graph obtained from $G^3$ by adding a new vertex $v$ adjacent to both degree $5$ vertices and one degree $3$ vertex. Thus $G$ contains a $W_4$. So $G$ is not $\langle\Z_{3}\rangle$-reduced by Lemma~\ref{neq0b}(i), a contradiction. Hence $|E(G)|\le 13$ and  $r(7, \Z_3)= 13$.~\qed

The proposition below  follows directly from the definitions of $r(n,\Z_3)$ and $\langle\Z_{3}\rangle$-reduced graphs.

\begin{proposition}\label{delta5andS}
  Let $G$ be a $\langle\Z_{3}\rangle$-reduced graph and $S\subset V(G)$. Then $$|\partial_G(S)|\ge \delta(G)|S| - 2 r(|S|, \Z_3).$$
\end{proposition}

By applying Lemma \ref{reducedgraph} and Proposition \ref{delta5andS} with straightforward calculation, we have the following lemma immediately.

\begin{lemma}
\label{5-edge-connected}
Let $G$ be a $\langle\Z_{3}\rangle$-reduced graph with $|V(G)| \leq 15$. If $\delta(G) \geq 5$, then $G$ is $5$-edge-connected and is essentially $8$-edge-connected. That is, for any $S\subset V(G)$ with $\min\{|S|, |\bar {S}|\}\ge 2,$ $$|\partial_G(S)|\ge 8.$$
\end{lemma}

We also need the following orientation theorem of Hakimi \cite{Haki65} to handle small graphs.

\begin{theorem} \label{Hakimi}{\em (Hakimi \cite{Haki65})}
Let $G$ be a graph and $\ell : V(G) \mapsto \Z$ be a function such that
$\sum_{v\in V(G)} \ell(v)=0$  and $\ell(v)\equiv d_G(v)\pmod 2, \; \forall v\in V(G)$.
Then the following are equivalent.
\\
(i) $G$ has an orientation $D$ such that $d_D^+(v)-d_D^-(v)=\ell(v)$, $\forall
v\in V(G)$.
\\
(ii)  $\displaystyle |\sum_{v\in S}\ell(v)|\le |\partial_G(S)|$, $\forall S \subset V(G)$.
\end{theorem}

\begin{lemma}\label{order13}
Every odd-$5$-edge-connected graph of order at most $13$ admits a mod $3$-orientation.
\end{lemma}

\proof ~Let $G$ be a counterexample with $|V(G)|+|E(G)|$ minimized. Then $G$ is a $\langle\Z_{3}\rangle$-reduced graph by Lemma \ref{kn}(iii). By Lemma \ref{splitting}, $G$ is $5$-regular, which implies that $|V(G)|$ is even. If $|V(G)|\le 10$, then $G$ has a mod $3$-orientation by Theorem \ref{fanzhou}, a contradiction. Assume $|V(G)|=12$ in the following.

 Since every even wheel is $\Z_3$-connected by Lemma \ref{neq0b}(i), $G$ does not contain an even wheel. If $G$ contains an odd wheel, then we apply $W$-contraction, and the resulting graph is still $5$-edge-connected by Lemma \ref{5-edge-connected} and Proposition \ref{PROP: Wcontraction}. This  yields a smaller counterexample by Lemma \ref{oddwheel}. Thus we obtain the following.

\medskip \noindent
{\bf Fact A.}  {\em $G$ does not contain a wheel as a subgraph. In particular, \\(i) $G$ contains no $K_4, W_4, W_5$;\\  (ii) for any vertex $v\in V(G)$, $G[N_G(v)]$ has no cycle, and therefore, $|E(G[N_G(v)])|\le 4$.}

\medskip

Let [X,Y] be a maximum edge cut of $G$  with $|X|\le |Y|$. Since $G$ is $5$-regular, we have
\begin{eqnarray} \label{eq1}
|[x,Y]_G|\ge 3,~\mbox{for any}~ x\in X~~\mbox{and}~~ |[y,X]_G|\ge 3,~\mbox{for any}~ y\in Y.
\end{eqnarray}
Hence
$$3|Y|\le |[X,Y]_G|\le 5|X|=5(12-|Y|),$$ which implies that $|Y|\le 7$ and thus $5\le |X|\le|Y|\le 7$ since $|X| + |Y| = 12$.

If $|X|=5$ and $|Y|=7$, denote $Y_0=\{y\in Y : |[y, X]_G|=3\}$. It follows that $5|X|-2|E(G[X])|=|[X,Y]_G|\ge 3|Y_0|+4(7-|Y_0|),$ which implies that
\begin{eqnarray}\label{XandY0}
  |Y_0|\ge 3+2|E(G[X])|~~\text{and}~~ |E(G[X])|\le 2.
\end{eqnarray}
By (\ref{XandY0}), there is a vertex $y_0 \in Y_0$ such that $y_0$ is adjacent to an isolated vertex in $G[X]$.  Since $y_0$ has only three neighbors in $X$, we have the following.

\medskip \noindent
 {\bf Fact B.}  {\em If $|X|=5$ and $|Y|=7$, then $N_G(y_0)\cap X$ induces a graph with at most one edge.}
\medskip

We define a function $\ell$ as follows.
If $|X|=|Y|=6$, set $\ell(x)=3$ for any $x\in X$ and $\ell(y)=-3$ for any  $y\in Y$; if $|X|=5$ and $|Y|=7$, set $\ell(x)=3$ for any  $x\in X$, $\ell(y_0)=3$ and $\ell(y)=-3$ for any $y\in Y\setminus\{y_0\}.$

As $\sum_{v\in V(G)} \ell(v)=0$ and by Theorem~\ref{Hakimi},  there exists an $S_0\subset V(G)$ with $|S_0|\le 6$ such that
\begin{eqnarray}\label{eqS3}
|\sum_{v\in S_0}\ell(v)|> |\partial_G(S_0)|.
\end{eqnarray}

Clearly, by (\ref{eq1}), we have
\begin{eqnarray}\label{notinXY}
 S_0\not\subseteq X~~ \text{and} ~~S_0\not\subseteq Y.
\end{eqnarray}

By (\ref{eqS3}) and Lemma \ref{5-edge-connected},  $|S_0|\ge 4$ and thus we have
 \begin{eqnarray}
 \label{eqS3L}
 |\sum_{v\in S_0}\ell(v)|>|\partial_G(S_0)|\ge 8.
 \end{eqnarray}

 We consider three cases according to $|S_0|$ in the following.

\noindent{\bf Case 1:} $|S_0|=4$.

Since $|\sum_{v\in S_0}\ell(v)|>|\partial_G(S_0)|\ge 8,$ by (\ref{eqS3L}), we have  $|\sum_{v\in S_0}\ell(v)|=12$. Thus, by (\ref{notinXY}), we have $|X|=5$, $|Y|=7$ and $S_0\cap Y=\{y_0\}$. However, it follows from Fact B  that
\begin{eqnarray}\label{contra}
 \partial_G(S_0)|\ge 5|S_0|-2|E(G[S_0])|\ge 20-2\times (1+3) = |\sum_{v\in S_0}\ell(v)|,
\end{eqnarray} a contradiction to (\ref{eqS3L}).

\noindent{\bf Case 2:} $|S_0|=5.$

With a similar calculation as in Case 1,  by Proposition \ref{delta5andS}, $|\partial_G(S_0)|\ge \delta(G)|S_0| - 2 r(5, \Z_3)\ge 9$, which implies $|\sum_{v\in S_0}\ell(v)|=15$ by (\ref{eqS3}). Thus $|X|=5$, $|Y|=7$ and $S\cap Y=\{y_0\}$ by (\ref{notinXY}). When $|X|=5$ and $|Y|=7$, we have $|E(G[X])|\le 2$ by (\ref{XandY0}). Since $|[y_0, X]_G|=3$ by the choice of $y_0$, we have  $|E(G[S_0])|\le 5$. Therefore, $\partial_G(S_0)|\ge 5|S_0|-2|E(G[S_0])|\ge 15$, which contradicts  (\ref{eqS3}).

\noindent{\bf Case 3:} $|S_0|=6.$

In this case, when $|S_0\cap X|= 2~\text{or}~3$, we have $|\sum_{v\in S_0}\ell(v)|\le 6$, a contradiction to (\ref{eqS3L}).
Thus $|S_0\cap X|= 1, 4, \text{or}~5$ by (\ref{notinXY}).

If $|S_0\cap X|=1~\text{or}~5$, then either $|S_0\cap X|=1$  or $|S_0\cap Y|  = 1$ and by  (\ref{eqS3}), we have  $$12\ge |\sum_{v\in S_0}\ell(v)|>|\partial_G(S_0)|=5|S_0|-2|E(G[S_0])|=30-2|E(G[S_0])|,$$
which implies $|E(G[S_0])|\ge 10.$

Let $w$ be the vertex in $G[S_0]$ such that $S_0\cap X=\{w\}$ or $S_0\cap Y=\{w\}$. Since $d_{G[S_0]}(w)\le d_{G}(w)=5$, we have $$|E(G[N_G(w)])|\ge |E(G[S_0])| - d_{G[S_0]}(w)\ge 5,$$
contradicting  Fact A(ii).

 If $|S_0\cap X|=4$, then $|S_0\cap Y|=2$. Since  $|\partial_G(S_0)|\ge 8$ by (\ref{eqS3L}), we have $|\sum_{v\in S_0}\ell(v)|= 12$,  implying that $|X|=5$, $|Y|=7$ and $y_0\in S_0\cap Y$. We claim that
 \begin{equation}\label{eq9}
   |E(G[S_0])|\le 9.
 \end{equation}

 If $|E(G[X])|= 2$, we have $|Y_0|=7$ by (\ref{XandY0}). Thus $|[S_0\cap X, S_0\cap Y]|\le 6$. Therefore $$|E(G[S_0])|\le |E(G[X])|+ |[S_0\cap X, S_0\cap Y]| +|E(G[S_0\cap Y])| \le 9.$$

 Now assume $|E(G[X])|\le 1$. Denote $(S_0\cap Y)\setminus\{y_0\}=\{z\}$. Since $|[y_0, X]_G|=3$ by the choice of $y_0$, we have $|[S_0\cap X, S_0\cap Y]| +|E(G[S_0\cap Y])|\le |[y_0, X]_G|+ d_G(z)\le 3+5$. Therefore, (\ref{eq9}) holds as well by the same inequality  above. This  proves (\ref{eq9}).

 By (\ref{eq9}), we have  $$12= |\sum_{v\in S_0}\ell(v)|>|\partial_G(S_0)|=5|S_0|-2|E(G[S_0])|\ge 30-18,$$
 a contradiction. This  contradiction completes the proof of the lemma. \qed

\begin{corollary}
\label{order15withk4}
  Let $G$ be a graph with $|V(G)|\le 15$. If $G$ is $5$-edge-connected and contains  a $K_4$, then $G$ admits a nowhere-zero $3$-flow.
\end{corollary}

\proof~  Let $G$ be a counterexample with $|E(G)|$ minimum. Denote $\{v_1, v_2, v_3, v_4\}$ to be the vertex set of a $K_4$ in $G$. We first show that $G$ is a $\langle\Z_{3}\rangle$-reduced graph. Suppose to the contrary that $H$ is a maximal nontrivial $\Z_3$-connected subgraph of $G$. Since $G/H$ admits no nowhere-zero $3$-flow, $G/H$ does not contain a $K_4$ by the minimality of $G$, and $|V(G/H)|\ge 14$ by Lemma \ref{order13}. So $|V(H)|=2$, meaning that $H$ consists of some parallel edges. Moreover one edge of $K_4$, say $v_1v_2$, is included in $H$. Then $V(H)=\{v_1, v_2\}$ and $H$ contains a digon $v_1v_2$.  This implies that $G[\{v_1, v_2, v_3, v_4\}]$ is $\Z_3$-connected by Lemma \ref{kn}(ii)(v), a contradiction to the maximality of $H$.
This proves that $G$ is a $\langle\Z_{3}\rangle$-reduced graph.

Since $G$ is $\langle\Z_{3}\rangle$-reduced and $\delta(G)\ge 5$,  $G$ is essentially $8$-edge-connected by Lemma \ref{5-edge-connected}. Applying $W$-contraction on $G$, by Proposition \ref{PROP: Wcontraction}, the resulting graph $G'$ remains $5$-edge-connected and has order at most $13$. By Lemma \ref{order13}, $G'$ admits a nowhere-zero $3$-flow. Therefore, $G$ admits a nowhere-zero $3$-flow by Lemma \ref{oddwheel}(i), a contradiction to the choice of $G$.
\qed

\section{Proofs of Theorems  \ref{mod3indle4odd}  and \ref{mod3indle3odd}}

This section will devote proofs of Theorems \ref{mod3indle4odd} and \ref{mod3indle3odd}. We start with some lemmas. For a vertex subset $X$ of $V(G)$, denote the neighbor set of $X$ to be  $N(X)=\{y | y\notin X~\mbox{and}~\mbox{there exists}~ x\in X ~\mbox{such that}~ xy\in E(G)\}$.

\begin{lemma}\label{Jlet-1}
  Let $G$ be a graph with $\alpha(G)\le t$.

  (i) For any nonempty $X\subset V(G)$,  $\alpha(G - (X\cup N(X)))\le t-1$.

  (ii) For any  maximal  $\Z_{3}$-connected subgraph $H$ of $G$ with $|V(H)|\ge t+1$,    $\alpha(G-V(H))\le t-1$.
\end{lemma}
\proof ~
  (i) is obvious.

  Now we prove (ii). Denote $J = G - V(H)$. Suppose to the contrary that $\alpha(J)= t$. Let $\{v_1, \dots, v_t\}$ be an independent set of size $t$ in $J$. By Lemma \ref{kn}(ii)(v), we have $|[v_i, V(H)]|\le 1$ for each $1\le i\le t$. Since $|V(H)|\ge t+1$, there exists a vertex $u\in V(H)$ such that $|[u, \{v_1, \dots, v_t\}]_G|=0$. Thus $\{v_1, \dots, v_t, u\}$ is an independent set of size $t+1$ in $G$, yielding a contradiction to $\alpha(G)\le t$.
\qed

\begin{lemma}\label{Ind3-14}
Let $G$ be a $\langle \Z_3\rangle$-reduced graph.

 (a) If $\alpha(G)\le 2$, then $|V(G)|\le 8$. Moreover, if $|V(G)|=8$, then $G$ contains  a $K_4$.

 (b) If $\alpha(G)\le 3$, then $|V(G)|\le 14$. Moreover, if $|V(G)|=14$, then $G$ is $5$-edge-connected and contains a $K_4$.

 (c) If $\alpha(G)\le 4$, then $|V(G)|\le 20$.
\end{lemma}

\proof ~  (a) Let $G$ be a $\langle\Z_{3}\rangle$-reduced graph  with $\alpha(G)\le 2$.  By Lemma \ref{kn}(vi), $G$ does not contain a $K_5$. Thus if $\kappa'(G) \in \{0, 1\}$, $|V(G)| \leq 8$ and if $|V(G)| = 8$, then $G$ contains a $K_4$.  (a) is also true by Theorem \ref{YangLi} if $\kappa'(G)\ge 3$.

Assume $\kappa'(G)=2$. Let $[X_1,X_2]_{G}$ be a 2-edge-cut of $G$, where $X_1,X_2$ form a  partition of $V(G)$. Then $|N(X_i)| \leq 2$ for each $i\in \{1,2\}$. By Lemmas \ref{kn}(vi) and \ref{Jlet-1}(i), for $i\in \{1, 2\}$, $G-(X_i\cup N(X_i))$ is a complete graph of size at most $4$. Since  $|N(X_i)|\le 2$, we have
$$|V(G)|\le 4+|X_i|+|N(X_i)|\le 6+|X_i|.$$
If $|X_i|\le 2$ for some $i\in \{1,2 \}$, then we have $|V(G)|\le 8$ and $|V(G)|= 8$ implies that $G$ contains a  $K_4$.
If both $|X_1|\ge 3$ and $|X_2|\ge 3$, then there is a vertex $x_1\in X_1$ such that $|[x_1, X_2]_{G}|=0$. Since $X_2\subseteq V(G)-(\{x_1\}\cup N(x_1))$,  by Lemma \ref{Jlet-1}(i),  $G[X_2]$ is a complete graph. Hence $|X_2|\le 4$ since $G$ does not contain a $K_5$. Similarly, we have $|X_1|\le 4$.
Thus $|V(G)|\le 8$. If  $|V(G)|= 8$, then $G[X_1]$ is a $K_4$.

(b)   By Lemma~\ref{min5}, $\delta(G) \leq 5$.  Thus by (a) and  Lemma~\ref{Jlet-1}(i) we have $|V(G)| \leq 1 + \delta(G)  + 8 \leq 14$ and if $|V(G)| = 14$, then $G$ contains a $K_4$. When $|V(G)|=14$, the above inequality is equality, implying every vertex in $G$ is of degree at least $5$. By Lemma \ref{5-edge-connected}, $G$ is $5$-edge-connected.

(c)  By Lemma~\ref{min5}, $\delta(G) \leq 5$.  Thus by (b) and  Lemma~\ref{Jlet-1}(i), we have $|V(G)| \leq 1 + \delta(G)  + 14  \leq 20$.
\qed

Now we are ready to prove  Theorem \ref{mod3indle3odd}.

\begin{customtheorem}{\ref{mod3indle3odd}}
  Every odd-$5$-edge-connected graph $G$  with $\alpha(G)\le 3$ admits a mod $3$-orientation.
\end{customtheorem}

 \medskip \noindent
{\bf Proof of Theorem \ref{mod3indle3odd}.} Let  $G$ be a counterexample with $|E(G)|$ minimum.   By Lemma \ref{splitting}, the degree of each vertex is odd; otherwise we lift all the edges incident with vertices of even degrees by applying Lemma \ref{splitting}, and then delete all isolated vertices to obtain a smaller counterexample. Thus $\delta(G)\geq 5$ and $G$ is $\langle\Z_{3}\rangle$-reduced.  By Lemma~\ref{order13}, $|V(G)|\ge 14$. Moreover, $|V(G)| \leq 14$ by Lemma \ref{Ind3-14}(b). Therefore $|V(G)| = 14$ and $G$ contains  a $K_4$. By Lemma~\ref{5-edge-connected},  $G$ is  $5$-edge-connected.
By Corollary \ref{order15withk4}, $G$ admits a mod $3$-orientation, a contradiction.
\qed

\begin{lemma}\label{reducetoorder15}
 Every  graph $G$ with $|V(G)|\ge 21$ and $\alpha(G)\le 4$ is $\langle\Z_{3}\rangle$-reduced to a graph of order at most $15$.
\end{lemma}
\proof~ Let $G_1$ be the underlying simple graph of $G$, which is obtained from $G$ by replacing parallel edges $[u,v]_G$  with a single edge $uv$ for each $|[u,v]_G|\ge 2$ in $G$. Since $|V(G_1)|\ge 21$,  $G_1$ contains a nontrivial $\Z_{3}$-connected subgraph by Lemma \ref{Ind3-14}(c), say $H_1$. Then $G[V(H_1)]$ is $\Z_{3}$-connected and $|V(H_1)|\ge 5$ by Lemma \ref{kn}(vi).  Let $H$ be a maximal  $\Z_{3}$-connected subgraph of $G$ containing $G[V(H_1)]$. Then we have $|V(H)|\ge |V(H_1)|\ge 5$.
 Let $J=G-V(H)$ and $J'$ be its $\langle\Z_{3}\rangle$-reduction. Since $|V(H)|\ge  5$, we have $\alpha(J)\le 3$ by Lemma \ref{Jlet-1}0(ii). Thus $J'$, the $\langle\Z_{3}\rangle$-reduction of $J$, is of order at most $14$ by Lemma \ref{Ind3-14}(b). Since  $|V(G')|=|V(J')|+1$, where $G'$ is the $\langle\Z_{3}\rangle$-reduction of $G$, we have $|V(G')|\le 15$. \qed

\medskip \noindent
{\bf  Equivalence of Theorems~\ref{mod3} and \ref{mod3indle4}:} In Theorem \ref{mod3}, if $|V(G)|\ge 21$, then the reduction of $G$ is of order at most $15$  by Lemma \ref{reducetoorder15}. So Theorem \ref{mod3indle4} is in fact equivalent to Theorem \ref{mod3} by Lemma \ref{order13}.

\medskip
Now we are ready to prove  Theorem \ref{mod3indle4odd}.
\begin{customtheorem}{\ref{mod3indle4odd}}
 Every odd-$5$-edge-connected graph $G$ of order at least $21$ with $\alpha(G)\le 4$ admits a mod $3$-orientation.
\end{customtheorem}

\medskip \noindent
{\bf Proof of Theorem \ref{mod3indle4odd}.} ~Let $G$ be a counterexample and $G'$ be its $\langle\Z_{3}\rangle$-reduction. We shall show that $G'$ has a mod $3$-orientation, which yields to a contradiction by Lemma \ref{kn}(iii).

 By Lemma \ref{reducetoorder15},  $|V(G')|\le 15$. Since $G'$ is odd-$5$-edge-connected with no mod $3$-orientation and by Lemma \ref{order13}, $|V(G')|\ge 14$. Therefore, $14\le |V(G')|\le 15$.

 Let $H$ be a maximal nontrivial $\Z_{3}$-connected subgraph of $G$ as in Lemma \ref{reducetoorder15}. Recall that $|V(H)|\ge 5$. Denote $v_1$ to be the contraction of $H$ in $G'$, and let $J'= G'- v_1$. Notice that  $J'$ is the $\langle\Z_{3}\rangle$-reduction of $J=G-V(H)$. Hence $\alpha(J')\le \alpha(J)\le 3$ by Lemma \ref{Jlet-1}(ii).

 We show the following to lead a contradiction.

 \medskip \noindent
 {\bf (I)} $|V(G')|= 14$.

  Suppose to the contrary $|V(G')|= 15$. Then $|V(J')|= 14$. So $J'$ is $5$-edge-connected and contains a $K_4$ by Lemma \ref{Ind3-14}(b).

   If $d_{G'}(v_1) < 5$, then $d_{G'}(v_1)$ is even since $G'$ is odd-$5$-edge-connected.  Applying Lemma \ref{splitting} to lift  all  edges incident with $v_1$, the resulting graph is $5$-edge-connected and of order $14$. By Corollary \ref{order15withk4}, the resulting graph admits a mod $3$-orientation and so does $G'$, a contradiction.

   If $d_{G'}(v_1)\geq 5$, then $G'$ is $5$-edge-connected as $J'$ is $5$-edge-connected. Since $G'$ contains a $K_4$, it admits a $3$-orientation by Corollary \ref{order15withk4}, a contradiction again. This proves (I).

\medskip \noindent
{\bf (II)}  $G'$ is $5$-regular and thus by Lemma~\ref{5-edge-connected}, $G'$ is $5$-edge-connected.

 Let $x$ be a vertex in $G'$. If $d(x)$ is even, applying Lemma \ref{splitting} to lift all the edges incident with $x$, the resulting graph remains odd-$5$-edge-connected with $13$ vertices. Thus it has a mod $3$-orientation by Lemma \ref{order13}, so does $G'$ by Lemma \ref{kn}(iii)(v), a contradiction. Thus $\delta(G') \geq 5$. By Lemma~\ref{5-edge-connected}, $G'$ is $5$-edge-connected  essentially $8$-edge-connected.

 Now assume $d_{G'}(x) \geq 7$. Since $\alpha(G') \leq 4$,  let $u,v$ be two adjacent vertices in $N_{G'}(x)$. Let $G''=G'_{[x, uv]}$ be the graph obtained from $G'$ by  deleting the edges $xu,xv$ and adding a new edge $uv$.   Since $G'$ is essentially $8$-edge-connected and $\delta(G'') \geq 5$, $G''$ remains $5$-edge-connected.  Note that  $G''$ contains a digon $uv$.  Then $G''/uv$ has $13$ vertices and remains $5$-edge-connected. Thus it has a mod $3$-orientation by Lemma \ref{order13}, so does $G'$ by Lemma \ref{kn}(iii)(v), a contradiction.

\medskip \noindent
{\bf The final step:}    By (II), $\delta(J') \leq 4$. Let $z\in V(J')$ with $d_{J'}(z) \leq 4$.  Since $\alpha(J') \leq 3$, by Lemma~\ref{Jlet-1}(i), $\alpha(J'-(\{z\}\cup N_{J'}(z))) \leq 2$.  Note that $J'-(\{z\}\cup N_{J'}(z))$ is a $\langle\Z_{3}\rangle$-reduced graph of order at least $8$. Thus by  Lemma~\ref{Ind3-14}(a), $J'-(\{z\}\cup N_{J'}(z))$ has exactly $8$ vertices and contains a $K_4$ and so does $G'$. By (I), (II) and Corollary \ref{order15withk4}, $G'$ admits a mod $3$-orientation, a contradiction. This completes the  proof of  Theorem \ref{mod3indle4odd}. \qed

\section*{Acknowledgement}
The second author was partially supported by  NSFC under grant numbers 11171288 and 11571149. The third author was partially supported by NSFC: 11401003.
We would like to thank Miaomiao Han and Professor Hong-Jian Lai for helpful discussion and valuable comments. We are grateful to the referees for their careful reading of the manuscript and helpful comments which led to the improvement of the presentation of this paper.

\end{document}